\documentclass{elsart}
\journal{Journal of Algebra}
\usepackage{amsmath,leftidx}
\usepackage{amsfonts}
\usepackage{amssymb}
\usepackage{amscd}

\def\diff{\ifmmode\Delta$-$\else$\Delta$-\fi}

\def\Q{{\mathcal Q}}

\def\Ob{{\mathcal Ob}}

\def\Z{{\mathbb Z}}

\def\C{{\mathbb C}}
\def\Q{{\mathbb Q}}

\DeclareMathOperator{\one}{\underline{\bf 1}}
\DeclareMathOperator{\GL}{\bf GL}
\DeclareMathOperator{\Gm}{\bf G_m}

\DeclareMathOperator{\Mn}{\bf M}
\DeclareMathOperator{\K}{\bf k}
\DeclareMathOperator{\Tr}{Tr}
\DeclareMathOperator{\A}{\mathcal A}

\DeclareMathOperator{\Hom}{Hom}
\DeclareMathOperator{\Span}{span}
\DeclareMathOperator{\End}{End}
\DeclareMathOperator{\Ind}{\mathcal Ind}
\DeclareMathOperator{\id}{id}
\DeclareMathOperator{\ev}{ev}
\DeclareMathOperator{\Sym}{Sym}

\DeclareMathOperator{\CoDiffT}{\bf CoDiff}

\DeclareMathOperator{\Aut}{\bf Aut}

\DeclareMathOperator{\Rep}{\bf Rep}
\DeclareMathOperator{\Vect}{\bf Vect}
\DeclareMathOperator{\Cat}{\mathcal{C}}

\DeclareMathOperator{\Char}{char}

\DeclareMathOperator{\AlgT}{{\bf Alg}_{\K}(\partial)}
\DeclareMathOperator{\Seq}{\mathcal{V}}

\newcommand{\Le}{\leqslant}
\newcommand{\Ge}{\geqslant}

\begin{document}

\begin{frontmatter}
\title{Differential Tannakian Categories\thanksref{thank1}}
\author{Alexey Ovchinnikov}
\ead{aiovchin@math.uic.edu}
\address{University of Illinois at Chicago\\ Department of Mathematics, Statistics, and Computer Science\\ Chicago, IL 60607-7045, USA}
\date\today
\thanks[thank1]{The work was partially supported by the Russian Foundation for Basic Research, project no. 05-01-00671 and by NSF Grant CCR-0096842.}


\begin{abstract}
We define a differential Tannakian category and show that under a
natural assumption it has
a fibre functor. If in addition this category is neutral, that is, the target category for
the fibre functor are finite dimensional vector spaces over the base field, then it is equivalent to the category of representations
of a (pro-)linear differential algebraic group. Our treatment of the problem is via differential Hopf algebras and Deligne's  fibre functor construction \cite{DeligneFS}. 
\end{abstract}

\end{frontmatter}

\section{Introduction}
A differential Tannakian category is a rigid abelian tensor category with an addition
differential structure. In \cite{DeligneFS} it is shown that if the dimension of each
object in a rigid abelian tensor category $\Cat$ over a field $\K$ of zero characteristic is a non-negative integer then there exists
a fibre functor $\omega$. This category is called neutral if the target category are finite dimensional
vector spaces over the field $\K$. In this case, the pair $(\Cat,\omega)$ is equivalent
to the category of finite dimensional representations of a (pro-)linear algebraic group \cite{DeligneFS,Saavedra,Deligne}. 

We find conditions for a rigid abelian tensor category $\Cat$ with an additional
differential structure over a 
field $\K$ of characteristic zero so that $\Cat$ has a fibre functor $\omega$ compatible with the differential structure. If, in addition,
this category is neutral then we show that $(\Cat,\omega)$ is equivalent to the category
of finite dimensional differential representations of a (pro-)linear differential algebraic
group.

Differential algebraic groups were introduced in \cite{Cassidy,CassidyRep} and further
used to build the differential Galois theory of parametrised linear differential equations \cite{PhyllisMichael}. This work was further generalised to include systems of linear difference equations with differential parameters~\cite{CharlotteMichael}.
Another approach to the Galois theory of systems
of linear differential equations with parameters is given in~\cite{TG}, where the authors
study Galois groups for generic values of the parameters. Also, it is shown in \cite{Hrushovsky} that over the field $\C(x)$  the differential Galois group of  a parametrised 
system of differential equations will be the same for all values of the parameter outside a countable 
union of proper algebraic sets.
 In the usual differential Galois theory the Tannakian approach is
a base for algorithms computing differential Galois groups which are linear algebraic groups \cite{Michael,Katz1,Katz2,Daniel1,Daniel2}. In the parametrised theory the Galois groups are linear {\it differential} algebraic
groups and the Tannakian theory that we develop should contribute to computation
of these groups.  

Tannakian approach for linear differential algebraic groups was first introduced in \cite{OvchRecoverGroup}, where a differential analogue of Tannaka's theorem
was proved (see~\cite{Springer,Water}). This was further developed in \cite{OvchTannakian}, where neutral
differential Tannakian categories were introduced. This approach used a given
fibre functor and the axioms for such a category were given heavily employing this
fibre functor. In \cite{OvchTannakian} it is shown that such a category is equivalent
to the category of finite dimensional differential representations of a (pro-)linear differential algebraic group and this is further applied to the Galois theory of parametrised linear differential  equations (see also~\cite{CharlotteMichael,Charlotte,CharlotteComp}). Generalisations of differential
Galois theory to non-linear equations via Galois groupoids can be found in
\cite{GuyCasale,GuyCasale2,GuyCasale3,GuyCasale4,GuyCasale5,Umemura,Umemura2,Malgrange,Malgrange2}.

In the present paper we extend the
results of \cite{OvchTannakian}. In particular, we no longer need the fibre functor to formulate
our axioms and, moreover, we show that under additional natural assumptions
such a functor exists. Note, that we do {\it not} require the base field to be differentially closed. The {\it main difficulty} was to give correct axioms relating
tensor product, duals, and the differential structure together without using the fibre functor and also to show
that the fibre functor that we construct ``respects'' the differential structure. 

A model theoretic treatment of general differential tensor categories 
and of some results of \cite{OvchTannakian} is given in~\cite{Moshe},
in particular, it is shown that a differential tensor category \cite[Definition 2.1]{Moshe}
together with a fibre functor compatible with the differential structure is equivalent
to the category of differential representations of the differential group scheme of
differential tensor automorphisms of the fibre functor. In \cite{Moshe}
the differential structure is given by a functor to the category of short exact sequences
and the tensor structure compatible with this derivation is introduced via Baer sums.

Our axioms are similar to the ones given in \cite{Moshe} but our version is more constructive. It is
assumed in \cite{Moshe} that the given fibre functor ``commutes'' with the
differential structure on the category, that is, for each object there is a functorial
``derivation'' satisfying the Leibniz rule with respect to the tensor structure. We do {\it not} make this  assumption on the fibre functor and
show how this functorial commutation for each object comes naturally from the category itself.

The paper is organised as follows. We review basic notions of differential algebra in Section~\ref{Basics}. In Section~\ref{MainDefinition} we introduce differential Tannakian
categories. In Remark~\ref{difftrivial} we show that the usual Tannakian categories can be given a trivial differential structure and, therefore, they are a part of our formalism. We recall what differential comodules and differential algebraic groups are in Section~\ref{secdiffcomod}. We show in Section~\ref{ConstructFiberFunctor}, Theorem~\ref{Th14}, how starting with a differential Tannakian category
one can get a fibre functor based on the ideas of \cite[Section 7]{DeligneFS}.
An additional assumption which guarantees that the category is neutral is made
in Section~\ref{additionalassumption}. We recover the differential Hopf algebra structure
of the corresponding group in Section~\ref{recoverbialgebra}. Lemmas~\ref{ComputationalDiffDual}~and~\ref{ProductRuleLemma} supply all necessary ingredients to finally recover the (pro-)linear differential algebraic group. In Proposition~\ref{LinearGroup} and Theorem~\ref{MainTheorem} of this section we  combine these lemmas with the results from~\cite{OvchRecoverGroup,OvchTannakian} to demonstrate the main reconstruction statement of the paper.

\section{Basic Definitions}\label{Basics}
A $\Delta$-ring $R$,
where $\Delta = \{\partial_1,\ldots,\partial_m\}$, is a commutative associative ring with unit $1$ and commuting derivations $\partial_i : R\to R$ such that
$$
\partial_i(a+b) = \partial_i(a)+\partial_i(b),\quad \partial_i(ab) =
\partial_i(a)b + a\partial_i(b)
$$
for all $a, b \in R$. If $\K$ is a field and
a $\Delta$-ring then $\K$ is called a $\Delta$-field. We restrict ourselves to the case of
$$
\Char\K = 0.
$$ If $\Delta = \{\partial\}$ then
we call a $\Delta$-field as $\partial$-field.
For example, $\Q$ is a $\partial$-field with the unique
possible derivation (which is the zero one). The field
$\C(t)$ is also a $\partial$-field with $\partial(t) = f,$
and this $f$ can be any rational function in $\C(t).$ Let $C$ be
the field of constants of $\K$, that is, $C = \ker\partial$.

Let 
$$
\Theta = \left\{\partial^i\:|\: i\in \Z_{\Ge 0}\right\}.
$$
Since $\partial$ acts on a $\partial$-ring $R$,
there is a natural action of $\Theta$ on $R$.
A non-commutative ring $R[\partial]$ of linear differential operators is generated as a left $R$-module by the monoid $\Theta$. A typical element
of $R[\partial]$ is a polynomial 
$$
D = \sum_{i=1}^na_i\partial^{i},\ a_i \in R.
$$
The right $R$-module structure follows from the
formula 
\begin{align}\label{basiccommutation}
\partial\cdot a = a\cdot \partial + \partial(a)
\end{align}
for all $a \in R$. Note that~\eqref{basiccommutation} defines both left and right
$R$-module structures on $R[\partial]$ that we will use in our constructions.
We denote the set of operators in $R[\partial]$ of order less than
or equal to $p$ by $R[\partial]_{\Le p}.$

Let $R$ be a $\partial$-ring. If $B$ is an $R$-algebra, then $B$ is a $\partial$-$R$-algebra
if the action of $\partial$ on $B$ extends the
action of $\partial$ on $R$. If $R_1$ and $R_2$
are $\partial$-rings then a ring homomorphism
$\varphi: R_1 \to R_2$ is called a $\partial$-homomorphism if it commutes with $\partial$, that is,
$$
\varphi\circ\partial = \partial\circ\varphi.
$$ We denote these homomorphisms simply by $\Hom(R_1,R_2)$.
If $A_1$ and $A_2$ are $\partial$-$\K$-algebras
then a $\partial$-$\K$-homomorphism simply means a 
$\K[\partial]$-homomorphism. We denote the
category of $\partial$-$\K$-algebras by $\AlgT$.
Let $Y = \{y_1,\ldots,y_n\}$ be a set of variables. We differentiate them:
$$
\Theta Y := \left\{\partial^iy_j
\:\big|\: i \in \mathbb{Z}_{\Ge 0},\ 1\Le j\Le n\right\}.
$$
The ring of differential polynomials $R\{Y\}$ in
differential indeterminates $Y$ 
over a $\partial$-ring $R$ is
the ring of commutative polynomials $R[\Theta Y]$
in infinitely many algebraically independent variables $\Theta Y$ with
the derivation $\partial$, which naturally
extends $\partial$-action on $R$ as follows:
$$
\partial\left(\partial^i y_j\right) := \partial^{i+1}y_j
$$
for all $i \in \Z_{\Ge 0}$ and $1 \Le j \Le n$.
A $\partial$-$\K$-algebra $A$ is called finitely
$\partial$-generated over $\K$ if there exists
a finite subset $X = \{x_1,\ldots,x_n\} \subset A$
such that $A$ is a $\K$-algebra generated by
$\Theta X$.  

\begin{defn}\label{SeqCategory}
The category $\Seq$ over a
$\partial$-field $\K$ is the category of
finite dimensional vector spaces over $\K$:
\begin{enumerate}
\item objects are finite dimensional $\K$-vector spaces, 
\item morphisms are $\K$-linear maps;
\end{enumerate}
with tensor product
$\otimes,$ direct sum $\oplus,$ dual $*,$ and an additional operation:
\begin{align}\label{seqdiff}
F : V \mapsto V^{(1)} := \leftidx{_{\K}}{\left(\left(\K[\partial]_{\Le 1}\right)_{\K}\otimes V\right)},
\end{align}
which we call a differentiation (or prolongation) functor. If $f \in \Hom(V,W)$ then we define
\begin{align}\label{seqmordiff}
F(f) : V^{(1)}\to W^{(1)},\ f(\partial^q\otimes v) = \partial^q\otimes f(v),\ 0\Le q\Le 1.
\end{align}
Here, $\K[\partial]_{\Le 1}$ is considered as the right $\K$-module of differential operators up to order $1$, $V$ is viewed as a left
$\K$-module, and the whole tensor product is viewed as a left $\K$-module as well.
\end{defn}

For each $V \in \Ob(\Seq)$ there are: a natural inclusion
\begin{align}\label{eqincl}
i_V: V \to V^{(1)},\quad v \mapsto 1\otimes v,
\end{align}
and a surjection
\begin{align}\label{eqdiff}
\varphi_V : V^{(1)} \to V,\quad v \mapsto 0,\: \partial\otimes v \mapsto v.
\end{align}
We also have a natural splitting of $\K^{(1)} = \K[\partial]_{\Le 1}\otimes \K$ as follows
\begin{align}\label{vonesplitting}
\K^{(1)} \ni a\otimes 1 + b\cdot\partial\otimes 1 \mapsto (a,b) \in \K\oplus\K
\end{align}
for all $a, b \in \K$. The projection onto the first component gives a left inverse to the
map $i_{\K}$. The other projection is just $\varphi_{\K}$.

\section{General definition of a differential Tannakian category}\label{DefinitionProperties}

Our general definition of differential Tannakian categories is motivated by the following example.
\begin{exmp} Let $K$ be a $\{\partial_t,\partial_x\}$-field with $\K$ as the $\partial_t$-field
of $\partial_x$-constants. Consider a system of ordinary differential equations with respect to the derivation $\partial_x$: 
\begin{equation}\label{pareq}
\partial_x Y = AY, 
\end{equation}
where $A \in \Mn_n(K)$ (see \cite{Michael}). Since $\partial_t$ and $\partial_x$ commute with each other,
system~\eqref{pareq} implies:
$$
\partial_x\partial_tY = \partial_t\partial_xY=\partial_t(AY) = (\partial_tA)Y+A(\partial_tY),
$$
which gives a new system ``prolonged'' with respect to $\partial_t$:
\begin{equation}\label{pareqprol}
\partial_x \begin{pmatrix}
\partial_t Y\\
Y 
\end{pmatrix}
= \begin{pmatrix}
A&\partial_t A\\
0&A
\end{pmatrix}
\begin{pmatrix}
\partial_tY\\
Y
\end{pmatrix}. 
\end{equation}
Equations~\eqref{pareq} and~\eqref{pareqprol} define two $K[\partial_x]$-modules \cite[Section 1.2]{Michael}
that we will denote by $M$ and $M^{(1)}$, respectively~\cite[Section 5]{OvchTannakian}. The matrix structure in~\eqref{pareqprol} implies
a canonical injective map $M \to M^{(1)}$ with a quotient $M^{(1)}/M \cong M$.
Thus, with a $K[\partial_x]$-module $M$ we associate the following exact sequence
of $K[\partial_x]$-modules:
\begin{equation*}
\begin{CD}
0@>>> M @>>>M^{(1)}@>>> M @>>> 0.
\end{CD}
\end{equation*}
As we will see, this is the main idea behind {\bf differential} Tannakian categories.
\end{exmp}

\subsection{Definition}\label{MainDefinition}
To define differential Tannakian categories  we need to recall a few constructions from
tensor categories. Let \begin{align}\label{commorphism}
\psi_{X,Y} : X\otimes Y \to Y\otimes X
\end{align}
be the {\it commutativity morphism} in $\Cat$ so that $\psi_{X,Y}\circ\psi_{Y,X} = \id_{Y,X}$ \cite[page 104]{Deligne}.
Let also
\begin{align}\label{DeltaFormula}
\Delta_X : X\otimes X^*\to X\otimes X^*\otimes X\otimes X^*
\end{align} be the dual morphism to the morphism that computes the evaluation \cite[page 110]{Deligne} of the leftmost $X^*$ with
the rightmost $X$ and switching the middle $X$ with $X^*$ that can be more formally written as 
\begin{align*}
( \ev_X\otimes\psi_{X,X^*})&\circ (\id_{X^*}\otimes\psi_{X,X}\otimes\id_{X^*})\circ(\id_{X^*\otimes X}\otimes\psi_{X^*,X}): X^*\otimes X\otimes X^*\otimes X\\
&\to X^*\otimes X\otimes X\otimes X^*\to X^*\otimes X\otimes X\otimes X^*\to X^*\otimes X
\end{align*}

\begin{defn}\label{NeutralParametricTannakian} A differential Tannakian category $\Cat$ over a 
field  $\K$ of characteristic zero is a
\begin{enumerate}
\item rigid
\item abelian
\item tensor
\end{enumerate}
category such that $\End(\one) = \K$, where $\one$ is the unit object with respect to the tensor product; together with:
\begin{enumerate}
\item a {\bf functor}
$$
F : X \mapsto X^{(1)},
$$  which 
is exact on the right and 
commutes with direct sums. 
\item
an inclusion (a morphism with the trivial kernel)
$$
i_X: X \to X^{(1)},
$$
which is a natural transformation from the $\id$ functor to the functor $F$;
\item
a morphism 
$$
\varphi_X : X^{(1)} \to X
$$ in $\Cat$, which is a natural transformation from the functor $F$  to the functor $\id$
such that
the short sequence 
\begin{equation}\label{diffinclusion}
\begin{CD}
0@>>>X @>{i_X}>> X^{(1)} @>{\varphi_X}>> X@>>> 0
\end{CD}
\end{equation}
is exact;
\item a chosen splitting 
\begin{equation}\label{onesplitting}
\begin{CD}
0@>>>\one @>{i_{\one}}>> \one\oplus\one @>{\varphi_{\one}}>> \one@>>> 0
\end{CD}
\end{equation}
\item {\bf functorial morphisms}: 
\begin{enumerate}
\item $$S_X : (X\otimes X^*)^{(1)} \to X^{(1)}\otimes \left(X^{(1)}\right)^*,$$
\item an injective (this is the Leibniz rule)
$$T_{X,Y}: (X\otimes Y)^{(1)} \to X^{(1)}\otimes Y^{(1)},
$$
\item and an isomorphism
$$D_X:  \left({X^{(1)}}\right)^* \to (X^*)^{(1)}$$
\end{enumerate} 
making diagrams~\eqref{morphismS}, \eqref{difftensordiagram2},  \eqref{difftensorcodiagram}, \eqref{difftensordiagram}, and \eqref{diffdual} given below commutative:
\end{enumerate}
\begin{equation}\label{morphismS}
\begin{CD}
 (X\otimes X^*)^{(1)}@>{S_X}>> X^{(1)}\otimes \left(X^{(1)}\right)^*\\
@V{F(\ev_X)}VV@V{\ev_{X^{(1)}}}VV\\
\one^{(1)}@>i_{\one}^{-1}>>\one
\end{CD}
\end{equation}
where $i_{\one}^{-1}$ denotes the left inverse to the map $i_{\one}$
in splitting~\eqref{onesplitting},
\begin{equation}\label{difftensordiagram2}
\begin{CD}
X\otimes Y @>{i_{X\otimes Y}}>>(X\otimes Y)^{(1)}@>{\quad\quad\quad\quad\varphi_{X\otimes Y}\quad\quad\quad\quad}>>X\otimes Y\\
@|@V{T_{X,Y}}VV@V(i_X\otimes\id_Y,\id_X\otimes i_Y)\circ\Delta VV\\
X\otimes Y@>{i_X\otimes i_Y}>>X^{(1)}\otimes Y^{(1)}@>{\left(\id_{X^{(1)}}\otimes\varphi_Y\right)\oplus\left(\varphi_X\otimes\id_{Y^{(1)}}\right)}>>X^{(1)}\otimes Y\oplus X\otimes Y^{(1)}
\end{CD}
\end{equation}
\begin{equation}\label{difftensorcodiagram}
\begin{CD}
(X\otimes X^*)^{(1)}@>{T_{X\otimes X^*,X\otimes X^*}\circ F(\Delta_X)}>>(X\otimes X^*)^{(1)}\otimes (X\otimes X^*)^{(1)}\\
@V{S_X}VV@V S_X\otimes S_X VV\\
X^{(1)}\otimes \left({X^{(1)}}\right)^*@>{\quad\quad\quad\Delta_{X^{(1)}}\quad\quad\quad}>>X^{(1)}\otimes \left({X^{(1)}}\right)^*\otimes X^{(1)}\otimes \left({X^{(1)}}\right)^*
\end{CD}
\end{equation}
\begin{equation}\label{difftensordiagram}
\begin{CD}
((X\otimes Y)\otimes(X\otimes Y)^*)^{(1)} @>{\alpha_{X,Y}\circ S_{X\otimes Y}}>>X^{(1)}\otimes Y^{(1)}\otimes\left((X\otimes Y)^{(1)}\right)^*\\
@VV{T_{X\otimes X^*,Y\otimes Y^*}\circ F\left(\id_X\otimes\psi_{Y,X^*}\otimes\id_{Y^*}\right)}V@AA\id_{X^{(1)}\otimes Y^{(1)}}\otimes{T_{X,Y}}^* A\\
(X\otimes X^*)^{(1)}\otimes(Y\otimes Y^*)^{(1)}@>{\beta_{X,Y}\circ\left(S_X\otimes S_Y\right)}>>X^{(1)}\otimes Y^{(1)}\otimes\left(X^{(1)}\otimes Y^{(1)}\right)^*\end{CD}
\end{equation}
where $$\alpha_{X,Y} := T_{X,Y}\otimes\id_{\left((X\otimes Y)^{(1)}\right)^*}$$
and
$$\beta_{X,Y} := \id_{X^{(1)}}\otimes\psi_{\left(X^{(1)}\right)^*,Y^{(1)}}\otimes\id_{\left(Y^{(1)}\right)^*};$$ 
\begin{equation}\label{diffdual}
\begin{CD}
\left(X\otimes X^*\right)^{(1)}@>\psi_{X^{(1)},\left(X^{(1)}\right)^*}\circ S_X>>\left(X^{(1)}\right)^*\otimes X^{(1)}\\
@VV S_{X^*}\circ F(\psi_{X,X^*})V@V D_X\otimes\id_{X^{(1)}}VV\\
\left(X^*\right)^{(1)}\otimes\left(\left(X^*\right)^{(1)}\right)^*@>\id_{\left(X^*\right)^{(1)}}\otimes {D_X}^*>>\left(X^*\right)^{(1)}\otimes X^{(1)}
\end{CD}
\end{equation}
\end{defn}

\begin{rem} Diagrams \eqref{morphismS}, \eqref{difftensordiagram2}, \eqref{difftensorcodiagram}, \eqref{difftensordiagram}, and \eqref{diffdual} are designed in such a way that the multiplication on the algebra of
regular functions that we recover satisfies the product rule and the coinverse is a differential
homomorphism. We will show this in Lemmas~\ref{ComputationalDiffDual} and~\ref{ProductRuleLemma}. 

Also, in Example~\ref{representationsexample} we will
explain why the category of differential representations (Definition~\ref{DiffComodules})
of a linear differential algebraic group (Definition~\ref{lineardefinition}) is a differential Tannakian category. In 
that example we will also argue why
the morphism $S_X$ is essential by showing that, in general, $$S_X \ne \left(\id_X\otimes (D_X)^{-1}\right)\circ T_{X,X^*}.$$
\end{rem}

\begin{rem} If $\Cat$ is $\Seq$ then exact sequence~\eqref{diffinclusion} splits. This
is in general not the case. In particular, if $\Cat$ is the category of differential representations of a linear differential algebraic group $G$ then sequence~\eqref{diffinclusion} splits (when $X$ is a faithful representation of $G$) if and only if the group $G$ is conjugate to a group
of matrices with constant entries \cite[Proposition 3]{OvchRecoverGroup}. The differential algebraic group $\Gm$ is not of that kind, for instance.
\end{rem}

\begin{rem}\label{difftrivial}More generally, the usual Tannakian categories do satisfy the formalism of
differential Tannakian categories given in Definition~\ref{NeutralParametricTannakian}.
Indeed, they can be supplied with the trivial differential structure:
$$
F : X \mapsto X\oplus X,
$$
$$
F(\varphi) := \varphi\oplus\varphi,\quad \varphi\in\Hom(X,Y),
$$
$$
S_X : (X\otimes X^*)\oplus(X\otimes X^*) \to (X\oplus X)\otimes(X^*\oplus X^*),
$$
$$
T_{X,Y} : (X\otimes Y)\oplus(X\otimes Y) \to (X\oplus X)\otimes(Y\oplus Y),
$$
and
$$
D_X : (X\oplus X)^*\to (X^*\oplus X^*)
$$
that all satisfy~ \eqref{morphismS}, \eqref{difftensordiagram2}, \eqref{difftensorcodiagram}, \eqref{difftensordiagram}, and \eqref{diffdual}, which
follows from the commutativity of $\oplus$ with $*$ and the distribution law
for $\oplus$ and $\otimes$.
\end{rem}

Similarly to~\cite[Section 2.2]{Moshe}:

\begin{lem}\label{diffstructlemma}
The field  $\K = \End(\one)$ has a natural {\it derivation} induced by the functor $F$.
\end{lem}
\begin{pf}
Let $a \in \K$. The functor $F$ defines a map $$\K \to \End\left(\one^{(1)}\right).$$ Since $\End\left(\one^{(1)}\right)$ is a $\K$-algebra, we have another map
$$G: \K \to \End\left(\one^{(1)}\right),\quad a \mapsto a\cdot\id_{\one^{(1)}}.$$ Since $i_{\one}$ is a natural transformation from $\id$ to $F$, the morphism $$(F-G)(a)$$ is zero on $i_{\one}(\one)$ for all $a \in \K$. Therefore, it induces a morphism $$H(a) : \one^{(1)}/i_{\one}(\one) \cong \one \to \one^{(1)}.$$ 
Moreover, $\varphi_{\one}\circ (F-G)(a)$ is the zero morphism  $\one^{(1)}\to \one$. Hence,
there exists a morphism $$D(a) : \one \to \one$$ such that $i_{\one}\circ D(a) = H(a)$.
Finally,
\cite[Claim, Section 2.2]{Moshe} shows that the  map  $a \mapsto D(a)$ is a {\it derivation} on $\K$.\qed 
\end{pf}

\begin{rem}Note that in the case of the trivial differential structure described in Remark~\ref{difftrivial} the derivation $D$ on $\K$ constructed above is the zero map as $F(a)=G(a)$ for all $a \in \K$ in this case.
\end{rem}

\begin{exmp} Let the prolongation functor $F$ act as in~\eqref{seqmordiff} on the morphisms. 
We will show how following Lemma~\ref{diffstructlemma} one can recover the differential
structure on the field $\K$ when, for example, $\Cat = \Rep_G$, where $G$ is a linear differential algebraic group (see Section~\ref{secdiffcomod}). So, for all $a,b \in \K$ we have
\begin{align*}
&F-G : \K \to \End(\K[\partial]),\quad a\mapsto (F-G)(a),\\
&(F-G)(a)(1\otimes b) = 1\otimes ab - a\otimes b = 0,\\
&(F-G)(a)(\partial\otimes b) = \partial\otimes ab - a\partial\otimes b = (\partial (a))\otimes b.
\end{align*}
Under the isomorphism $\one^{(1)}/i_{\one}(\one)\cong \one$ we have $\partial\otimes b \mapsto b$. This composed with $i_{\one}$ sends $\partial\otimes b \mapsto 1\otimes b$. Then, for each $a \in \K$ the functor $F-G$ induces the map
$$
H(a) : \one \to \one^{(1)},\quad b\mapsto (\partial(a))\otimes b, 
$$
which induces the map
$$
D(a) : \one \to \one,\quad b \mapsto (\partial(a))\cdot b,
$$
that is, as an element of $\K$, $D(a) =\partial(a)$.
\end{exmp}

\subsection{Differential comodules}\label{secdiffcomod}
 Let $A$ be a $\partial$-$\K$-algebra.
Assume that $A$ is supplied with the following operations:
\begin{itemize}
\item differential algebra homomorphism $m : A\otimes A \to A$ is the
multiplication map on $A,$
\item differential algebra homomorphism $\Delta : A \to A\otimes A$
which is a comultiplication,
\item differential algebra homomorphism $\varepsilon : A  \to \K$
which is a counit,
\item differential algebra homomorphism $S : A \to A$
which is a coinverse.
\end{itemize}
We also assume that these maps satisfy commutative diagrams (see \cite[page 225]{CassidyRep}):
\begin{equation}\label{HopfAxioms}
\begin{CD}
A @>{\Delta}>> A\otimes A\\
@V{\Delta}VV @V{\id_A\otimes\Delta}VV \\
A\otimes A @>{\Delta\otimes \id_A}>> A\otimes A\otimes A
\end{CD}\qquad
\begin{CD}
A @>{\Delta}>> A\otimes A\\
@V{\id_B}VV @V{\id_A\otimes\varepsilon}VV \\
A @>{\sim}>> A\otimes \K
\end{CD}\quad\quad
\begin{CD}
A @>{\Delta}>> A\otimes A\\
@V{\varepsilon}VV @V{m\circ(S\otimes \id_A)}VV \\
\K @>{\hookrightarrow}>> A
\end{CD}
\end{equation}

\begin{defn}\label{DiffHopfAlgebra} Such a commutative associative $\partial$-$\K$-algebra $A$
with unit $1$ and operations $m,$ $\Delta,$ $S,$
and $\varepsilon$ satisfying axioms~\eqref{HopfAxioms} is called a {\it differential Hopf algebra} (or $\partial$-$\K$-Hopf algebra).
\end{defn}

\begin{defn}\label{DiffComodules} A  finite dimensional vector space $V$ over $\K$ is called an {\it $A$-comodule}
if there is a given $\K$-linear morphism
$$
\rho : V \to V\otimes A,
$$
satisfying the axioms:
\begin{equation*}
\begin{CD}
V@>{\rho}>> V\otimes A\\
@VV{\rho}V @VV{\id_V\otimes\Delta}V \\
V\otimes A@>{\rho\otimes\id_A}>> V\otimes A\otimes A
\end{CD}\qquad\qquad
\begin{CD}
V @>{\rho}>> V\otimes A\\
@VV{\id_V}V @VV{\id_V\otimes\varepsilon}V \\
V @>{\sim}>> V\otimes \K
\end{CD}
\end{equation*}
If $A$ is a differential coalgebra with $\Delta$ and $\varepsilon$
then for a comodule $V$ over $A$ the $\K$-space $\K[\partial]_{\Le 1}\otimes V = V^{(1)}$ has a natural $A$-comodule structure
\begin{align}\label{DiffFla}
\rho(f\otimes v) := f\otimes\rho(v)
\end{align}
for  $f \in \K[\partial]_{\Le 1}$ and $v \in V.$
\end{defn}

We denote the category of  comodules over a differential coalgebra
$A$ by $\CoDiffT_A$ with the induced differentiation~\eqref{DiffFla}.
For a $\partial$-$\K$-Hopf algebra $A$ the functor
$$
G : \AlgT \to \{\mathrm{Groups}\}, \quad R \mapsto \Hom(A,R)
$$
is called an {\it affine differential algebraic group scheme} generated by $A$ (see \cite[Section 3.3]{OvchRecoverGroup}). In this case $V \in \Ob(\CoDiffT_A)$ is called a
{\it differential representation} of $G$ (\cite[Definition 7, Theorem 1]{OvchRecoverGroup}).
The category $\CoDiffT_A$ is also denoted by $\Rep_G$.

The differential $\GL_n$ by definition is the functor 
represented by the $\partial$-$\K$-Hopf algebra
$$
\K\{X_{11},\ldots,X_{nn},1/\det(X)\},
$$
where $X_{ij}$ are differential indeterminates.
The comultiplication $\Delta$ and coinverse  $S$ are defined on $X_{ij}$ in the usual way. Their
prolongation on the derivatives of $X_{ij}$ can be
obtained by differentiation.

\begin{defn}\label{lineardefinition} An affine differential algebraic group scheme $G$ is called a {\it linear
differential algebraic group} if there exists
an embedding:
$$
G \to \GL_n
$$
for some $n \in \Z_{\Ge 1}$.
\end{defn}

\begin{exmp}\label{representationsexample} Let $G$ be a linear
differential algebraic group. The category $\Rep_G$ is a rigid abelian tensor category
with $\End(\one)=\K$ and $\K$ is a differential field. The differential structure on $\Rep_G$
is given by~\eqref{seqdiff},~\eqref{seqmordiff}, and~\eqref{DiffFla}. Therefore, the $\K$-linear
maps~\eqref{eqincl} and~\eqref{eqdiff} are $G$-morphisms. Splitting~\eqref{onesplitting} comes from splitting~\eqref{vonesplitting}. It remains to show that there exist $G$-equivariant maps 
$$S_X, T_{X,Y},\text{and}\ \: D_X$$ for all $X,Y \in\Ob\left(\Rep_G\right)$.
The morphism $S_X$ is defined in~\cite[Lemma 4, part (3)]{OvchRecoverGroup}.
For $T_{X,Y}$ and $D_X$ this is shown in~\cite[Lemmas 7 and 11]{OvchRecoverGroup},
respectively. We will now show that, in general,
$$S_X \ne \left(\id_X\otimes (D_X)^{-1}\right)\circ T_{X,X^*}.$$
Indeed, in the notation of Section~\ref{recoverbialgebra} and \cite{OvchRecoverGroup} for $X = \Span\{v_1,\ldots,v_n\}$ and $X^*=\Span\left\{v^1,\ldots,v^n\right\}$ with $v^j(v_i) = \delta_{i,j}$, $1\Le i,j\Le n$, we have
\begin{align*}
\ev_X\left(S_X\left(\partial\otimes\left(v_i\otimes v^j\right)\right)\right) &= \ev_X\left(\left(\partial\otimes v_i\right)\otimes F\left(v^j\right)\right) = F\left(v^j\right)(\partial\otimes v_i) = \\&=\partial\left(v^j(v_i)\right) =
 \partial(\delta_{i,j}) = 0.
\end{align*}
where
$$
F: X^* \to \left(X^{(1)}\right)^*,\: F(u)(v) = u(v),\: F(u)(\partial\otimes v) = \partial(u(v)),\: v \in X,\: u\in X^*.
$$
On the other hand (see \cite[proof of Lemma 11]{OvchRecoverGroup} for $(D_X)^{-1}$),
\begin{align*}
&\ev_X\circ\left(\id_X\otimes (D_X)^{-1}\right)\circ T_{X,X^*}\left(\partial\otimes\left(v_i\otimes v^j\right)\right) =\\
&= \ev_X\circ\left(\id_X\otimes (D_X)^{-1}\right)\left(\partial\otimes v_i\otimes 1\otimes v^j+1\otimes v_i\otimes \partial\otimes v^j\right) =\\
&=\ev_X\left(\partial\otimes v_i\otimes (\partial\otimes v_j)^*+1\otimes v_i\otimes F\left(v^j\right)\right) = \delta_{i,j}+\delta_{i,j} \ne 0.
\end{align*}
\end{exmp}

\section{Main result}\label{mainresultsection}
\subsection{Constructing a differential fibre functor}\label{ConstructFiberFunctor} Let $\Cat$ be a differential Tannakian category over a $\partial$-field $\K$ of characteristic zero. As in \cite[Theorem 7.1]{DeligneFS}, {\bf suppose} that for each $X\in\Ob(\Cat)$ we have 
$$
\dim X \in \Z_{\Ge 0},
$$
where $\dim X = \Tr(\id_X)$ and $\Tr : \End(X) = \Hom(\one,X^*\otimes X) \to \End(\one)$ given by the evaluation morphism (see \cite[formula (1.7.3)]{Deligne}).
As it is shown in \cite[Section 7.18]{DeligneFS}, one can construct a ring object $A$ of the
category $\Ind\Cat$ (see \cite[page 167]{DeligneFS}) with $\mu : A\otimes A\to A$ and $e : \one\to A$ such that for all $X \in \Ob(\Cat)$
we have
\begin{align}\label{isomtoA}
A\otimes X \cong A^{\dim X}
\end{align}
and 
\begin{align}\label{fibrefunctor}
\omega(X) := \Hom(\one,A\otimes X)
\end{align}
is a fibre functor. However, this functor does not necessarily preserve the differential
structure, that is, one does not necessarily have a functorial isomorphism $$\omega\left(X^{(1)}\right) \to \omega(X)^{(1)}.$$ Therefore, in order to achieve our goal we have to modify the functor. Based on $A$, we shall construct another ring object $B$ of $\Ind\Cat$ giving the required functorial isomorphism.
\begin{thm}\label{Th14}
There exists a ring object $B$ of $\Ind\Cat$ such that if we define the functor $\omega$  by 
\begin{equation}\label{fibrefunctorB}
\omega(X) = \Hom(\one,B\otimes X),
\end{equation} $X\in\Ob(\Cat)$, we will have:
\begin{enumerate}
\item $\omega$ is a fibre functor over $\Hom(\one,B)$ in the sense of \cite[Section 1.9]{DeligneFS},
\item there exists a functorial isomorphism $\omega\left(X^{(1)}\right)\to \omega(X)^{(1)}$, $X \in \Ob(\Cat)$.
\end{enumerate}
\end{thm}
\begin{pf} 
Define a ``differential polynomial'' ring object (cf. \cite[Section 1.2]{Henri})
$$
B := \bigcup_{n\Ge 1,\:p\Ge 0 }\Sym^n\left(F^p(A)\right).
$$
To continue with our argument we need the following
\begin{lem} We have:
\begin{enumerate}
\item
there exists  an isomorphism
\begin{equation}\label{functorialisom}
d_B : \omega\left(B^{(1)}\right) \cong \omega(B)^{(1)}
\end{equation}
functorial with respect to $\End(B)$,
\item $B$ is faithfully flat over $A$,
\item tensoring with $B$ splits short exact sequences from $\Cat$,
\item for each $X \in \Ob(\Cat)$ there exists an isomorphism
\begin{equation}\label{isomtoB}
B\otimes X \cong B^{\dim X}.
\end{equation}
\end{enumerate}
\end{lem}
\begin{pf}
There is a functorial isomorphism
\begin{equation}\label{FuncIsom}
(B\otimes B)^{(1)} \to B\otimes B^{(1)}
\end{equation}
given by iterated applications of functorial injections 
$T_{\cdot,\cdot}$
defined in~\eqref{difftensordiagram2}. Indeed, 
diagram~\eqref{difftensordiagram2} implies by induction an injection
for each $n$ and $p$ 
$$
(\Sym^n(F^p(A)))^{(1)} \to \Sym^n(F^{p+1}(A)),
$$
making the following diagram commutative:
\begin{equation*}
\begin{CD}
(\Sym^n(F^p(A)))^{(1)} @>>> \Sym^n(F^{p+1}(A))\\
@AAA@AAA\\
\Sym^n(F^p(A))@=\Sym^n(F^p(A))
\end{CD}
\end{equation*}
which after taking the union with respect to all $n$ and $p$ induces
a functorial isomorphism 
\begin{equation}\label{BtoB}
B^{(1)} \to B.
\end{equation} A similar argument provides
a functorial isomorphism $$(B\otimes B)^{(1)}\to B\otimes B,$$
which combined with~\eqref{BtoB} gives  isomorphism~\eqref{FuncIsom}.
Now, there is also a functorial isomorphism
$$
\Hom(\one,Y)^{(1)} \to \Hom\left(\one,Y^{(1)}\right)
$$
for any object $Y$.
Setting $Y = B\otimes B$, this gives us~\eqref{functorialisom}. 
Let 
$$
0\to X\to Y\to Z \to 0
$$
be a short exact sequence in $\Cat$. By \cite[Section 7.18]{DeligneFS}, the sequence
$$
0 \to A\otimes X\to A\otimes Y\to A\otimes Z\to 0
$$
splits. Then, the functors 
\begin{equation}\label{functors}
\left(F^p\right)^{\otimes n}\otimes\id_{\Cat},\quad p\Ge 0,\: n\Ge 1,
\end{equation} 
give a splitting of the exact sequence
$$
0 \to B\otimes X\to B\otimes Y\to B\otimes Z\to 0.
$$
Using functors~\eqref{functors} and additivity of $F$, isomorphism~\eqref{isomtoB} follows from~\eqref{isomtoA}. 
As in the proof of \cite[Lemma 7.15]{DeligneFS},
the ring object $B$ is faithfully flat over $A$, and, therefore, we have constructed $B$ in $\Ind(\Cat)$ with the desired properties.
\qed
\end{pf}
Statement (2--4) of the lemma imply statement (1) of the theorem (see \cite[Sections 1.9 and 7.18]{DeligneFS}).
It follows from~\eqref{isomtoB} that the functorial in $X$ morphism of $B$-modules 
\begin{align}\label{surjectivetoAX}
\pi_1 : B^{(\oplus\Hom(B,B\otimes X))}\to B\otimes X
\end{align}
is surjective. Here, for a set $I$, $B^{(\oplus I)}$ denotes the direct sum of $B$ with
respect to $I$.
The functorial isomorphism chosen in~\eqref{functorialisom} 
induces (by taking direct sums) a functorial isomorphism  
\begin{equation}\label{functorialA}
d_{B^{(\oplus I)}}  : \omega\left(\left(B^{(\oplus I)}\right)^{(1)}\right) \to \omega\left(B^{(\oplus I)}\right)^{(1)}.
\end{equation}
By taking a functorial in $X$ free resolution of~\eqref{surjectivetoAX}: 
\begin{equation}\label{freeresolution}
\begin{CD}
B^{(\oplus I_X)}@>{\pi_2}>> B^{(\oplus\Hom(B,B\otimes X))}@>{\pi_1}>> B\otimes X@>>> 0,
\end{CD}
\end{equation}
where $I_X := \Hom\left(B,\ker\left(B^{(\oplus\Hom(B,B\otimes X))}\to B\otimes X\right)\right)$,
one gets a functorial (in $X$) isomorphism
\begin{equation}\label{AAXsequence}
d_X :  \omega\left(X^{(1)}\right) \to \omega(X)^{(1)}.
\end{equation}
Indeed, the functor $F$ is exact on the right and, therefore, transforms free resolution~\eqref{freeresolution}
of $B\otimes X$ to a free resolution of $(B\otimes X)^{(1)}$ as follows:
\begin{equation*}
\begin{CD}
 B^{(\oplus I_X)}@>>> \left(B^{(\oplus I_X)}\right)^{(1)}@>>>B^{(\oplus I_X)}\\
@V{\pi_2}VV@V F(\pi_2)VV@V{\pi_2}VV\\
B^{(\oplus\Hom(B,B\otimes X))} @>{}>> \left(B^{(\oplus\Hom(B,B\otimes X))}\right)^{(1)} @>{}>> B^{(\oplus\Hom(B,B\otimes X))}\\
@V{\pi_1}VV@V{F(\pi_1)}VV@V{\pi_1}VV\\
B\otimes X @>>> (B\otimes X)^{(1)} @>>> B\otimes X\\
@VVV@VVV@VVV\\
0@.0@.0
\end{CD}
\end{equation*}
Applying the exact functor $\omega$ and using the isomorphism given in~\eqref{functorialA},
we have:
\begin{equation*}
\begin{CD}
\omega\left(\left(B^{(\oplus I_X)}\right)^{(1)}\right)@>d_{B^{(\oplus I_X)}}>>\omega\left(B^{(\oplus I_X)}\right)^{(1)}\\
@VVV@VVV\\
\omega\left(\left(B^{(\oplus\Hom(B,B\otimes X))}\right)^{(1)}\right) @>d_{B^{(\oplus\Hom(B,B\otimes X))}}>> \omega\left(B^{(\oplus\Hom(B,B\otimes X))}\right)^{(1)}\\
@VVV@VVV\\
\omega\left((B\otimes X)^{(1)}\right) @. \omega(B\otimes X)^{(1)}\\
@VVV@VVV\\
0@.0
\end{CD}
\end{equation*}
where the upper square is commutative because the isomorphism $d_B$ is functorial.
Therefore, one can define a functorial isomorphism 
\begin{equation}\label{derivationAX}
d_{B\otimes X} :  \omega\left((B\otimes X)^{(1)}\right) \to \omega(B\otimes X)^{(1)}
\end{equation}
By recalling~\eqref{fibrefunctor} and restricting in~\eqref{derivationAX} to $$\mathcal{I}_X : X \to \one\otimes X \to B\otimes X,$$ employing the morphism 
$$\mu\otimes\id_X : B\otimes B\otimes X \to B\otimes X$$ 
with $\left(\mu\otimes\id_X\right)\circ \left(\id_B\otimes \mathcal{I}_X\right)= \id_{B\otimes X}$, the isomorphism $d_{B\otimes X}$ induces
 a functorial in $X$ isomorphism $d_X$ in~\eqref{AAXsequence}. Indeed, for
 $$\varphi \in \omega\left(X^{(1)}\right) = \Hom\left(\one,B\otimes X^{(1)}\right)$$ we have 
\begin{align*}
\omega&(F(\id_B\otimes\mathcal{I}_X))\circ\varphi \in \Hom\left(\one,B\otimes (B\otimes X)^{(1)}\right) = \omega\left((B\otimes X)^{(1)}\right) \to\\
&\to\omega(B\otimes X)^{(1)} 
=\Hom(\one,B\otimes B\otimes X)^{(1)}\xrightarrow{\circ\omega\left(F\left(\mu\otimes\id_X\right)\right)}\\
&\to\Hom(\one,B\otimes X)^{(1)} =\omega(X)^{(1)}.
\end{align*}
\qed
\end{pf}

\subsection{Additional assumption}\label{additionalassumption}
Assume now that
\begin{align}\label{assumption}
\Hom(\one,B)\cong \K.
\end{align}
It then follows from isomorphism~\eqref{isomtoB} that $\omega(X)$ is a finite dimensional vector space over $\K$ for
all $X \in \Ob(\Cat)$ making $\omega$ a fibre functor over $\K$; and 
we have a  functorial isomorphism $$d_X : \omega\left(X^{(1)}\right) \cong \omega(X)^{(1)},$$
such that the following diagram is commutative:
\begin{equation*}
\begin{CD}
\omega(X) @>{\omega\left(i_{X}\right)}>> \omega\left(X^{(1)}\right) @>{\omega\left(\varphi_{X}\right)}>> \omega(X)\\
@|@V d_X VV @ |\\
\omega(X) @>{i_{\omega(X)}}>>\omega(X)^{(1)} @>{\varphi_{\omega(X)}}>> \omega(X)
\end{CD}
\end{equation*}

\subsection{Equivalence of categories}\label{recoverbialgebra} 
Recall (from \cite[Definition 8]{OvchRecoverGroup}) that the functor $\Aut^{\otimes,\partial}(\omega)$ from $\partial$-$\K$-algebras to groups is defined as follows. For a $\partial$-$\K$-algebra $K$ the group $\Aut^{\otimes,\partial}(\omega)(K)$ is the set of sequences
$$\lambda(K) = (\lambda_X\:|\: X\in \Ob(\Cat)) \in
\Aut^{\otimes,\partial}(\omega)(K)$$ such that $\lambda_X$ is a $K$-linear automorphism
of $\omega(X)\otimes K$ for each object $X$, that is,
$\lambda_X \in \Aut_K(\omega(X)\otimes K)$, such that
\begin{itemize}
\item for all $X_1,$ $X_2$ we have
\begin{align}\label{TensorSpreading}
\lambda_{X_1\otimes X_2} = \lambda_{X_1}\otimes\lambda_{X_2},
\end{align}
\item $\lambda_{\underline{1}}$ is the identity map on
$\omega(\underline{1})\otimes K$,
\item for every $\alpha \in \Hom(X,Y)$ we have
\begin{align}\label{Equivariance}
\lambda_Y\circ(\alpha\otimes \id_K) = (\alpha\otimes \id_K)\circ\lambda_X :
\omega(X)\otimes K \to \omega(Y)\otimes K,
\end{align}
\item for every $X$ we have
\begin{align}\label{CommuteWithD}
F(\lambda_X) = \lambda_{X^{(1)}},
\end{align}
\item the group operation $\lambda_1(K)\cdot\lambda_2(K)$ is defined
by composition in each\\ $\Aut_K(\omega(X)\otimes K)$.
\end{itemize}

\begin{thm}\label{MainTheorem}
Let $\Cat$ be a differential Tannakian category and $\omega : \Cat \to \Vect_{\K}$
be the fibre functor defined by~\eqref{fibrefunctorB} with additional assupmtion~\eqref{assumption}.
Then
$$
(\Cat,\omega) \cong \Rep_G
$$
for the differential group scheme
$$G = \Aut^{\otimes,\partial}(\omega),$$
which is a pro-linear differential algebraic group, that is,
$G$ as a functor is represented by a direct limit of finitely $\partial$-generated 
$\partial$-$\K$-Hopf
algebras.
\end{thm}
\begin{pf} The proof will consist of several steps,  lasting until the end of the paper finishing with
Proposition~\ref{LinearGroup}, which combined with \cite[Corollary 2]{OvchTannakian} gives the result.

First, we will
recall the constructions given in \cite[Section 3.4]{OvchTannakian}.
Let $X$ be an object of $\Cat$ and $\{\{X\}\}$ (respectively, $\Cat_X$) be the full abelian (respectively, full abelian tensor) subcategory of $\Cat$ generated by $X$
(respectively, containing $X$ and closed under the functor $F$).
Consider
$$
F_X = \bigoplus_{V \in \Ob(\{\{X\}\})}  V\otimes\omega(V)^*.
$$
This tensor product is understood in the sense of \cite[page 131]{Deligne}. For an object $V$ of the category $\{\{X\}\}$ we have the canonical injections:
$$
j_V : V\otimes\omega(V)^* \to F_X.
$$
Consider the minimal subobject $R_X$ of $F_X$ with subobjects
\begin{align*} \left\{(j_V(\id\otimes \phi^*) -
j_W(\phi\otimes \id))(V\otimes
\omega(W)^*)\:\big|\:V,W \in\Ob(\{\{X\}\}),\phi\in
\Hom(V,W)\right\}.
\end{align*}
We let
$$
P_X = F_X\big/R_X.
$$
We now put
$$
A_X = \omega(P_X),
$$
which is a $\K$-vector space.
Let $V$ be an object of $\{\{X\}\}.$
For
$$
v  \in \omega(V),\
u \in \omega(V)^*
$$
we denote by
$$
a_V\left(v\otimes u\right)
$$
the image in $A_X$ of
$$
\omega(j_V)\left(v\otimes u\right).
$$
So, for any $\phi \in \Hom(V,W)$ we have
\begin{align}\label{morphismsanda}
a_V\left(v\otimes \omega(\phi)^*(u)\right) =
a_W(\omega(\phi)(v)\otimes u).
\end{align}
for all
$$
v \in \omega(V),\ u \in \omega(W)^*.
$$

Let us define a {\bf comultiplication} on $A_X$.
Let $\{v_i\}$ be a basis of $\omega(V)$ with the dual basis
$\{u_j\}$ of $\omega(V)^*.$
We let
\begin{equation}\label{Comult}
\Delta : a_V(v\otimes u)
\mapsto\sum_i a_V(v_i\otimes u)\otimes a_V(v\otimes u_i),
\end{equation}
giving the same as~\eqref{DeltaFormula}.
The {\bf counit} is defined in the following way:
\begin{equation}\label{Counit}
\varepsilon : a_V(v\otimes u) \mapsto \omega(\ev_{V})(v\otimes u) = u(v).
\end{equation}
The {\bf coinverse} is defined as follows. We let
\begin{equation}\label{Coinverse}
S : a_V(v\otimes u) \mapsto a_{V^*}(\omega(\psi_{V,V^*})(v\otimes u)) = a_{V^*}(u\otimes v). 
\end{equation}

We shall introduce a {\bf differential algebra} structure on 
$$\A_X = \varinjlim\limits_{Y \in \Ob(\Cat_X)} A_Y$$
as follows. Let $V, W \in \Ob(\Cat_X)$, $v \in \omega(V)$, $w \in \omega(W)$, $u \in \omega(V)^*$, and $t \in \omega(W)^*$. 
Define a {\bf multiplication} on $\A_X$ by
\begin{align}\label{multiplicationformula}
a_V(v\otimes u)\cdot a_W(w\otimes t) := a_{V\otimes W}((v\otimes w)\otimes(u\otimes t)).
\end{align} 

Before we define a derivation on $\A_X$ we will do some preparation.
Since  $T_{\:\cdot\:,\:\cdot}$ is functorial in each argument, the following diagram
is commutative:
\begin{equation}\label{Tandone}
\begin{CD}
 ((V\otimes V^*)\otimes (W\otimes W^*))^{(1)}@>{T_{V\otimes V^*,W\otimes W^*}}>> (V\otimes V^*)^{(1)}\otimes (W\otimes W^*)^{(1)}\\
@VV{F(\ev_V\otimes\ev_W)}V@V{F(\ev_V)\otimes F(\ev_W)}VV\\
\left(\one\otimes\one\right)^{(1)}@>T_{\one,\one}>>\one^{(1)}\otimes\one^{(1)}
\end{CD}
\end{equation}
It follows from~\eqref{difftensordiagram2} that
\begin{align}\label{omegaT}
\omega\left(T_{V\otimes V^*,W\otimes W^*}\right)&(\partial\otimes(v\otimes u\otimes w\otimes t))= \partial\otimes(v\otimes u)\otimes1\otimes (w\otimes t)+\\
&+1\otimes (v\otimes u)\otimes\partial\otimes(w\otimes t) + \notag\\ 
&+1\otimes(\Phi_{V\otimes V^*}\otimes\id+\id\otimes\Phi_{W\otimes W^*})((v\otimes u)\otimes(w\otimes t)),\notag
\end{align}
where $\Phi$ is a functorial $\K$-liner map.

\begin{rem} According to~\eqref{omegaT} the morphism $T_{V\otimes V^*,W\otimes W^*}$ gives us
the product rule for the tensor structure up to the linear term. We will take this linear term
given by $\Phi$ into account defining the differential structure on $\A_X$ in~\eqref{defdiffstructure}.
\end{rem}

\begin{lem} We have
\begin{align}\label{phieval}
\omega(\ev_V)(\Phi_{V\otimes V^*}(v\otimes u))=0.
\end{align}
\end{lem}
\begin{pf}
Indeed, diagram~\eqref{difftensordiagram2} with $X=Y=\one$ together with splitting~\eqref{onesplitting} give us
$$
\omega(T_{\one,\one})(\partial\otimes(u(v)\otimes t(w)))= \partial\otimes u(v)\otimes \omega(i_{\one})(t(w))+\omega(i_{\one})(u(v))\otimes\partial\otimes t(w).$$
Therefore,
$$
\omega(\ev_V\otimes \ev_W)\left((\Phi_{V\otimes V^*}\otimes\id+\id\otimes\Phi_{W\otimes W^*})((v\otimes u)\otimes(w\otimes t))\right) = 0.
$$
Since $\Phi$ is functorial, taking $W := V$, $w := v$, $t := u$, we obtain the result.\qed
\end{pf}

We define a {\bf derivation} on $\A_X$ as follows:
\begin{align}\label{defdiffstructure}
\partial(a_V(v\otimes u)) := a_{V^{(1)}}\left(\omega(S_V)(\partial\otimes(v\otimes u)+2\otimes\Phi_{V\otimes V^*}(v\otimes u))\right).
\end{align}

\begin{lem}\label{ComputationalDiffDual} 
Let $V\in \Ob(\Cat_X)$, $v \in \omega(V)$, and $u\in\omega(V)^*$. Then, we have
\begin{align}\label{diffdualcommuteformula}
&a_{{(V^*)}^{(1)}}\left(\omega\left(S_{V^*}\circ F\left(\psi_{V,V^*}\right)\right)(\partial\otimes(v\otimes u)+2\otimes\Phi_{V\otimes V^*}(v\otimes u))\right) = \\
&=a_{{\left(V^{(1)}\right)}^*}\left(\omega\left(\psi_{V^{(1)},\left(V^{(1)}\right)^*}\circ S_V\right)\left(\partial\otimes(v\otimes u)+2\otimes\Phi_{V\otimes V^*}(v\otimes u)\right)\right)\notag
\end{align}
and the ``differential evaluation'' is
\begin{align}\label{diffeval}
\omega\left(\ev_{V^{(1)}}\circ S_V\right)(\partial\otimes(v\otimes u)+2\otimes\Phi_{V\otimes V^*}(v\otimes u)) = \partial(u(v)).
\end{align}
\end{lem}
\begin{pf}
To show~\eqref{diffdualcommuteformula} note that using~\eqref{morphismsanda} we obtain that
\begin{align*}
&a_{{(V^*)}^{(1)}}\left(\omega\left(S_{V^*}\circ F\left(\psi_{V,V^*}\right)\right)(\partial\otimes(v\otimes u)+2\otimes\Phi_{V\otimes V^*}(v\otimes u))\right) = \\
&=a_{(V^*)^{(1)}}(\omega(S_{V^*})(\partial\otimes(u\otimes v)+2\otimes\Phi_{V\otimes V^*}(v\otimes u))) =\\
&=  a_{\left(V^{(1)}\right)^*}\left(\omega\left(\left(D_V^{-1}\otimes {D_V}^*\right)\circ S_{V^*}\right)(\partial\otimes(u\otimes v)+2\otimes\Phi_{V\otimes V^*}(v\otimes u))\right),
\end{align*} 
since $D_V$ is an isomorphism. Moreover, from~\eqref{diffdual} we conclude that
\begin{align*}
&a_{{\left(V^{(1)}\right)}^*}\left(\omega\left(\psi_{V^{(1)},\left(V^{(1)}\right)^*}\circ S_V\right)(\partial\otimes(v\otimes u)+2\otimes\Phi_{V\otimes V^*}(v\otimes u))\right) = \\
&= a_{\left(V^{(1)}\right)^*}\left(\omega\left(\left(D_V^{-1}\otimes {D_V}^*\right)\circ S_{V^*}\right)(\partial\otimes(u\otimes v)+2\otimes\Phi_{V\otimes V^*}(v\otimes u))\right),
\end{align*}
which proves~\eqref{diffdualcommuteformula}.
Let us show~\eqref{diffeval} now. According to~\eqref{morphismS} and~\eqref{phieval} we have 
\begin{align*}
\omega&\left(\ev_{V^{(1)}}\circ S_V\right)(\partial\otimes(v\otimes u)+2\otimes\Phi_{V\otimes V^*}(v\otimes u)) =\\
&= \omega\left(i_{\one}^{-1}\circ F(\ev_V)\right)\left(\partial\otimes(v\otimes u)+2\otimes\Phi_{V\otimes V^*}(v\otimes u)\right) = \\
&= \omega\left(i_{\one}^{-1}\right)\left(\partial\otimes u(v)+0\right) = \omega\left(i_{\one}^{-1}\right)\left(\partial(u(v))\otimes 1 + u(v)\cdot\partial\otimes 1\right) =\\
&= \partial(u(v))
\end{align*}
because $\omega\left(i_{\one}^{-1}\right)$ is the projection to the first copy of $\K$ in
splitting~\eqref{vonesplitting}.\qed
\end{pf}

\begin{rem}\label{structureremark} Now, it follows from~\cite[Lemma 11]{OvchRecoverGroup} that our coinverse $S$ is a $\partial$-$\K$-algebra
homomorphism. 
Formula~\eqref{diffeval} shows that the counit $\varepsilon$
is a $\partial$-$\K$-algebra homomorphism as well. It remains to
show that the comultiplication $\Delta$ is a $\partial$-$\K$-algebra homomorphism and that the differential structure we defined satisfies the product rule.  We will do this in Lemmas~\ref{CoproductLemma} and~\ref{ProductRuleLemma}. Correctness of $\Delta$ is
shown in~\cite[Lemma 9]{OvchRecoverGroup}. The Hopf algebra commutative diagrams
for $\Delta$, $S$, and $\varepsilon$ are justified in~\cite[Lemmas 9--11]{OvchRecoverGroup}.
\end{rem}

\begin{lem}\label{CoproductLemma} The comultiplication $\Delta$ defined in~\eqref{Comult} is a 
$\partial$-$\K$-algebra homomorphism of $\A_X$.
\end{lem}
\begin{pf}
Follows from diagram~\eqref{difftensorcodiagram} and the observations given in Remark~\ref{structureremark}.\qed
\end{pf}

\begin{lem}\label{ProductRuleLemma} Formula~\eqref{defdiffstructure} defines a differential algebra structure on $\A_X$, that is,
$$
\partial\left(a_V(v\otimes u)\cdot a_W(w\otimes t)\right) = \partial(a_V(v\otimes u))\cdot a_W(w\otimes t) + a_V(v\otimes u)\cdot\partial(a_W(w\otimes t))
$$
for all $V, W \in \Ob(\Cat)$, $v\in\omega(V)$, $u \in\omega(V)^*$, $w \in \omega(W)$,
and $t \in\omega(W)^*$. 
\end{lem}
\begin{pf}
By definition,
\begin{align*}
&\partial \left(a_V(v\otimes u)\cdot
a_W(w\otimes t)\right) = \partial \left(a_{V\otimes W}((v\otimes w)\otimes(u\otimes t))\right) =\\
&= a_{(V\otimes W)^{(1)}}
(\omega(S_{V\otimes W})(\partial\otimes((v\otimes w)\otimes(u\otimes t))+\\
&\quad\quad+2\otimes\Phi_{(V\otimes W)\otimes(V\otimes W)^*}((v\otimes w)\otimes(u\otimes t)))).
\end{align*}
We have 
\begin{align*}
&\left(\partial a_V(v\otimes u)\right)\cdot
a_W(w\otimes t) + a_V(v\otimes u)\cdot
\left(\partial a_W(w\otimes t)\right)=\\
&=a_{V^{(1)}}(\omega(S_V)(\partial\otimes(v\otimes u)+2\otimes\Phi_{V\otimes V^*}(v\otimes u)))\cdot a_W(w\otimes
t)+\\
&\quad+a_V(v\otimes u)\cdot
a_{W^{(1)}}(\omega(S_W)(\partial\otimes(w\otimes t)+2\otimes\Phi_{W\otimes W^*}(w\otimes t))) =\\
&=a_{V^{(1)}\otimes W^{(1)}}(\omega(S_V\otimes S_W)((\partial\otimes(v\otimes u)+\\&\quad\quad\quad\quad\quad+2\otimes\Phi_{V\otimes V^*}(v\otimes u))\otimes \omega(i_W)(w\otimes t)+\\
&\quad\quad\quad\quad\quad+\omega(i_V)(v\otimes u)\otimes(\partial\otimes(w\otimes t)+2\otimes\Phi_{W\otimes W^*}(w\otimes t))))=\\
&=a_{V^{(1)}\otimes W^{(1)}}(\omega((S_V\otimes S_W)\circ T_{V\otimes V^*,W\otimes W^*})(\partial\otimes((v\otimes u)\otimes(w\otimes t))+\\
&\quad\quad\quad\quad\quad+2\otimes\Phi_{(V\otimes W)\otimes(V\otimes W)^*}((v\otimes w)\otimes(u\otimes t)))).
\end{align*}
The last equality follows from~\eqref{omegaT} 
and the facts that $\Phi$ is functorial and the category $\Cat$ has an associativity morphism for
the tensor product:
$$
\Phi_{V\otimes V^*}\otimes\id_{W\otimes W^*} + \id_{V\otimes V^*}\otimes\Phi_{W\otimes W^*} = \left(\id_V\otimes\psi_{W,V^*}\otimes\id_{W^*}\right)\circ\Phi_{(V\otimes W)\otimes(V\otimes W)^*}.
$$
Now, the product rule we are showing follows from diagram~\eqref{difftensordiagram}, formula~\eqref{morphismsanda},
and the condition that $T_{V,W}$ is injective.\qed
\end{pf}

We have all ingredients by now to
recover a (pro-)linear differential algebraic group from $(\Cat,\omega)$.
\begin{prop}\label{LinearGroup} The group $G_X$, defined by $G_X(R) = \Hom(\A_X,R)$ for each $\partial$-$\K$-algebra $R$,
is a linear differential algebraic group.
\end{prop}
\begin{pf}
It follows from Lemmas~\ref{ComputationalDiffDual},~\ref{CoproductLemma}, and ~\ref{ProductRuleLemma} and Remark~\ref{structureremark} combined with \cite[Proposition 2]{OvchTannakian} and \cite[Corollary 1]{OvchTannakian} that $\A_X$ is a
finitely $\partial$-generated $\partial$-$\K$-Hopf algebra. In particular,  our Lemmas~\ref{ComputationalDiffDual} and \ref{ProductRuleLemma} are used to substitute 
\cite[Lemma 6]{OvchTannakian} and \cite[Lemmas 7 and 11]{OvchRecoverGroup} employed in \cite[Proposition 2]{OvchTannakian}.\qed
\end{pf}
Then, $G = \varprojlim G_X$ is represented by $\varinjlim \A_X$ and this finishes the proof of Theorem~\ref{MainTheorem}.\qed
\end{pf}

\section{Acknowledgements} The author is highly grateful to Ben Antieau, Bojko Bakalov, Daniel Bertrand, Pierre Deligne, Henri Gillet, Sergey Gorchinskiy,
Christian Haesemeyer, Charlotte Hardouin, Moshe Kamensky, Claudine Mitschi, Jacques Sauloy, Michael Singer, Dima Trushin, and participants of the Kolchin Seminar for their helpful comments
and support. The author also thanks the referees for their important suggestions.

\bibliographystyle{elsart-num}
\bibliography{difftanncat}

\end{document}